# The V/S test of long-range dependence in random fields

## Frédéric Lavancier

*Université de Nantes, Laboratoire Jean Leray*
*UFR Sciences et Techniques*
*2 rue de la Houssinière - BP 92208 -*
*F-44322 Nantes Cedex 3, France*
*e-mail:* frederic.lavancier@univ-nantes.fr

**Abstract:** Recently, Giraitis et al. (2003, [10]) proposed the $V/S$ statistic for testing long memory in random sequences. We generalize this statistic to the setting of random fields. The null hypothesis is concerned with short memory random fields while the alternative contains a very large family of long memory random fields. Contrary to most of the previous works dealing with long-range dependence, no assumption is made about the isotropy of the strong dependence. Some simulations are presented in order to assess the power of the test according to the kind of long memory in presence.

## 1. Introduction

A stationary random field $X = (X_n)_{n \in \mathbb{Z}^d}$ is usually said to exhibit long memory when its covariance function $r(n)$, $n \in \mathbb{Z}^d$, is not absolutely summable: $\sum_{n \in \mathbb{Z}^d} |r(n)| = \infty$. An alternative definition relates on spectral properties: A random field is said to exhibit long memory if its spectral density is unbounded. These two points of view are closely related but not equivalent. In this paper, the concepts of "strong dependence" and "long-range dependence" are the same as "long memory".

Most of the previous studies on long memory random fields (see [6, 7, 18]) assume that the strong dependence occurs with the same intensity in all directions. Indeed, these works are concerned with isotropic long memory according to the following definition.

**Definition 1.** A stationary random field exhibits isotropic long memory if it admits a spectral density which is continuous everywhere except at 0 where

$$f(x) \sim ||x||^{\alpha - d} b\left(\frac{x}{||x||}\right) L\left(\frac{1}{||x||}\right), \quad 0 < \alpha < d,$$

where $||.||$ denotes the Euclidean norm, where $L$ is slowly varying at infinity and $b$ is a continuous function on the unit sphere in $\mathbb{R}^d$.

However, it is easy to construct non-isotropic long memory random fields. In [14], such fields arise from particular filterings of a white noise, from the aggregation of weakly dependent random fields, or from systems of statistical mechanics in phase transition.

Our aim is to construct a procedure for discriminating between weak dependent random fields and strong dependent ones, regardless of the isotropy.

In dimension $d = 1$, several tests for long memory are available. They are mostly based on an estimation of the variations of the partial sums process of $X$. For all these tests, the alternative hypothesis consists in parametric families of long memory processes, typically FARIMA time series. Lo (1991, [17]) first developed a test based on the $R/S$ statistic, which estimates the range of the partial sums process of $X$. The KPSS test was initially developed by Kwiatkovski et al. (1992, [12]) for testing stationarity (under weak dependence assumptions) against the presence of a trend or a unit root. A variant of the KPSS test, based on an estimation of the second order moments of the partial sums of $X$, was proposed by Lee and Schmidt (1996, [16]) in order to test long memory. From the same idea, Giraitis et al. (2003, [10]) introduced the $V/S$ statistic, based on an estimation of the variance of the partial sums of $X$. It appears that this test is more powerful than the $R/S$ test and than the $KPSS$ test for detecting strong dependence.

Note that the behavior of the partial sums of $X$ is not the unique way to test long memory in dimension $d = 1$. Goodness-of-fit tests for long range dependent time series have been developed (cf. [1, 9] and [4]). To the best of our knowledge, in dimension $d > 1$, there exists no generic model able to exhibit different situations such as isotropic or non-isotropic long memory. A goodness-of-fit approach seems therefore too restrictive in our framework.

In this paper, we generalize the $V/S$ test to the setting of random fields. Since this test is based on the partial sums of $X$, we summarize in Section 2 what is known about their limiting behavior under short memory and long memory. This will allow us to specify our testing hypothesis. Section 3 proves the consistency of the test through asymptotic results. Some simulations are presented in Section 4 when $d = 2$. They reveal that the power of the test is strongly related to the anisotropy of the long memory.

## 2. Test hypotheses

Let $X$ be a second order stationary real random field. We want to test the null hypothesis: *X is weakly dependent* against the alternative assumption: *X exhibits long memory*. The V/S test that we extend in this paper does not rely exactly on these testing hypotheses. It rather focuses on the behavior of the partial sums process, defined for all $t \in [0,1]^d$ by

$$\sum_{k \in A_n(t)} X_k,$$

with $A_n(t) = \mathbb{Z}^d \cap \prod_{i=1}^{d}[1, \lfloor (n-1)t_i \rfloor + 1]$, where $\lfloor x \rfloor$ denotes the integral value of $x$ and $n$ is a positive integer. Indeed, there is a close relation between the dependency of $X$ and the asymptotic behavior of its partial sums.

When $X$ is a centered and weakly dependent random field, i.e. when its covariance function is absolutely summable, it is well known that if $\sigma^2 := \sum_{h \in \mathbb{Z}^d} r(h) \neq 0$, a functional central limit theorem generally holds:

$$\frac{1}{\sigma n^{d/2}} \sum_{k \in A_n(t)} X_k \stackrel{\mathcal{D}([0,1]^d)}{\longrightarrow} B(t). \qquad (2.1)$$

Here, $B$ denotes the Brownian Sheet, i.e. the centered Gaussian process such that $E(B(t)B(s)) = \prod_{i=1}^{d} t_i \wedge s_i$. The convergence takes place in $\mathcal{D}([0,1]^d)$, the Skorokhod space of cadlag functions defined on $[0,1]^d$ (see [2] for instance). This result has been proved by different authors according to the kind of weak dependence of $X$, among others: Wichura (1969, [20]) when $X$ is an i.i.d process, Dedecker (2001, [5]) under a weak projective assumption.

On the other hand, when $X$ is long-range dependent, (2.1) is generally false. More precisely, when the spectral density of $X$ exhibits at least one singularity located at zero, then, under some structural hypotheses on $X$, it has been proved that

$$\frac{1}{n^{\gamma}} \sum_{k \in A_n(t)} X_k \stackrel{\mathcal{D}([0,1]^d)}{\longrightarrow} Y(t), \qquad (2.2)$$

where $\gamma > d/2$ and $Y$ is a random field different from the Brownian Sheet (not even necessarily Gaussian). This result is shown in [6] for functionals of Gaussian fields, extended in [18] to functionals of linear fields, when the long memory is isotropic according to Definition 1. In the case of linear fields which exhibit anisotropic long memory, (2.2) is proved in [15].

Let us now precise the testing hypotheses.

**H0: Short memory hypothesis.** The second order random field $X$ is stationary, with a covariance function $r$, such that **SM1**, **SM2** and **SM3** below are satisfied.

**SM1**

$$\sum_{j \in \mathbb{Z}^d} |r(j)| < \infty \quad \text{and} \quad \sigma^2 := \sum_{j \in \mathbb{Z}^d} r(j) > 0. \qquad (2.3)$$

**SM2**

$$\frac{1}{\sigma n^{d/2}} \sum_{k \in A_n(t)} (X_k - E(X_0)) \stackrel{\mathcal{D}([0,1]^d)}{\longrightarrow} B(t),$$

where $B$ is the Brownian Sheet.

**SM3**

$$\sup_{i \in \mathbb{Z}^d} \sum_{(j,k) \in \mathbb{Z}^{2d}} |c_4(i,j,k)| < \infty,$$

where $c_4$ represents the fourth order cumulants of $X$: Denoting $\tilde{X}_i = X_i - E(X_i)$, $c_4(i,j,k) = E[\tilde{X}_0 \tilde{X}_i \tilde{X}_j \tilde{X}_k] - r(i)r(k-j) - r(j)r(k-i) - r(k)r(j-i)$.

**H1: Long memory hypothesis.** The second order random field $X$ is stationary and satisfies **LM1** and **LM2** below.

**LM1**

$$\frac{1}{n^\gamma L(n)} \sum_{k \in A_n(t)} (X_k - E(X_0)) \xrightarrow{\mathcal{D}([0,1]^d)} Y(t), \tag{2.4}$$

with $\gamma > d/2$, where $L$ is a slowly varying function at infinity and $Y$ is some measurable and non-degenerated random field.

**LM2**

$$Var\left(\sum_{k \in A_n(1)} X_k\right) = O(n^{2\gamma} L^2(n)). \tag{2.5}$$

*Remark* 1. The null hypothesis **H0** deals with short memory as suggested by assumption **SM1**. Assumption **SM2** claims that a functional central limit theorem holds, as expected in the short memory setting. However, **SM1** does not imply **SM2**. For instance, some counter-examples are available in Theorems 7 and 8 in [11]. Finally, **SM3** is needed for technical reasons.

As explained above, (2.4) holds for a large family of long memory fields and it can be reasonably chosen as the alternative hypothesis. However, when the strong dependence involves non zero spectral singularities, then (2.1) may remain true (cf. [15], Theorem 2 and Theorem 4). Then, there is no chance for this particular situation of long memory to be detected by the test. Notice that the same restriction exists in dimension 1 in all the long memory testing procedures based on the behavior of the partial sums. Assumption **LM2** is convenient for technical reasons and appears to be a weak restriction to **LM1**.

## 3. The testing procedure

### 3.1. Test statistic

We generalize the V/S statistic to $d > 1$. Let us first introduce some notations. For all positive integer $n$, let $A_n \equiv A_n(1)$. We denote, for all positive integer $j$,

$$S_{n,j}^* = \sum_{i \in A_j} \left(X_i - \overline{X}_n\right), \tag{3.1}$$

where $\overline{X}_n = n^{-d} \sum_{j \in A_n} X_j$. Let $q$ be an integer in $[1, n]$. An estimator of $\sigma^2$, defined by (2.3), is

$$\hat{s}_n^2 = \sum_{j \in B_{q-1}} \omega_{q,j} \hat{r}(j), \tag{3.2}$$

with $B_q = \{-q, \ldots, q\}^d$. Here $\omega_{q,j} = \prod_{i=1}^d (1 - \frac{|j_i|}{q})$ are some weights leading to the positivity of $\hat{s}_n^2$ (see for instance [3] p360) and $\hat{r}$ is the empirical covariance function:

$$\hat{r}(j) = \frac{1}{n^d} \sum_{k_1=1}^{n-|j_1|} \cdots \sum_{k_d=1}^{n-|j_d|} \left(X_{k_1,\ldots,k_d} - \overline{X}_n\right) \left(X_{k_1+|j_1|,\ldots,k_d+|j_d|} - \overline{X}_n\right).$$

The statistic V/S is defined by

$$M_n = n^{-d} \frac{\widehat{Var}\left(S_{n,j}^*, \ j \in A_n\right)}{\hat{s}_n^2},$$

where $\widehat{Var}\left(S_{n,j}^*, \ j \in A_n\right) = n^{-d} \sum_{j \in A_n} \left(S_{n,j}^* - \overline{S_n^*}\right)^2$ and $\overline{S_n^*} = n^{-d} \sum_{j \in A_n} S_{n,j}^*$. One can rewrite $M_n$ as

$$M_n = \frac{n^{-2d}}{\hat{s}_n^2} \left[ \sum_{j \in A_n} S_{n,j}^{*}{}^2 - n^{-d} \left( \sum_{j \in A_n} S_{n,j}^* \right)^2 \right]. \tag{3.3}$$

### *3.2. Consistency of the test*

In what follows, the integer $q$ involved in Definition 3.2 is actually a function $q_n$ depending on $n$. The following proposition establishes the consistency of the test when the sequence $q_n$ is properly chosen.

**Proposition 1.** *If* $\lim_{n \to \infty} q_n = \infty$ *and* $\lim_{n \to \infty} q_n / n = 0$, *then,*

*(i) Under* **H0**,

$$M_n \xrightarrow{\mathcal{L}} \int_{[0,1]^d} \left( B(t) - \left( \prod_{i=1}^d t_i \right) B(1) \right)^2 dt - \left[ \int_{[0,1]^d} \left( B(t) - \left( \prod_{i=1}^d t_i \right) B(1) \right) dt \right]^2, \tag{3.4}$$

*where $B$ is the Brownian Sheet on $[0,1]^d$ and $\xrightarrow{\mathcal{L}}$ denotes the convergence in law.*

*(ii) Under* **H1**,

$$M_n \xrightarrow{P} \infty, \tag{3.5}$$

*where $\xrightarrow{P}$ denotes the convergence in probability.*

The testing procedure is the following: Given a significance level $\alpha \in [0,1]$, one rejects the null hypothesis **H0** if $M_n$, given by (3.3), is greater than $c(\alpha)$, where $c(\alpha)$ is such that

$$P(U_d > c(\alpha)) = \alpha,$$

where $U_d$ is distributed according to the asymptotic law involved in (3.4). Proposition 1 guarantees that the significance level of the test is asymptotically correct and that the power of the test goes to 1.

*Remark* 2. To our knowledge, the theoretical form of the asymptotic law $U_d$ in (3.4) is unknown for $d > 1$. When $d = 1$, the distribution function of $U_1$ is $F_K(\pi\sqrt{x})$, where $F_K$ denotes the Kolmogorov distribution function (see [10]). This identification comes from Watson (1961, [19]) and is not easy to extend to $d > 1$. It is however straightforward to obtain $E(U_d) = (1/2)^d + (1/4)^d - 2(1/3)^d$. For $d = 2$, a simulation of the density distribution of $U_2$ is presented in Section 4.1 (cf. Figure 1).

*Proof.* Notice that if $\frac{k}{n} \leq t_1 < \frac{k+1}{n}$, then

$$S_n^*(\lfloor nt_1 \rfloor + 1, \ldots, \lfloor nt_d \rfloor + 1) = S_n^*(k+1, \lfloor nt_2 \rfloor + 1, \ldots, \lfloor nt_d \rfloor + 1),$$

where $S_n^*(j) \equiv S_{n,j}^*$ is defined for all $j \in A_n$ in (3.1). Therefore

$$n^{-d} \sum_{j \in A_n} S_{n,j}^* = \int_{[0,1]^d} S_n^*(\lfloor nt_1 \rfloor + 1, \ldots, \lfloor nt_d \rfloor + 1) dt.$$

The same equality holds with respect to $S_n^{*2}$ and expression (3.3) of $M_n$ becomes

$$M_n = \frac{n^{-d}}{\hat{s}_n^2} \left[ \int_{[0,1]^d} S_n^{*2}(\lfloor nt_1 \rfloor + 1, \ldots, \lfloor nt_d \rfloor + 1) dt \right.$$
$$\left. - \left( \int_{[0,1]^d} S_n^*(\lfloor nt_1 \rfloor + 1, \ldots, \lfloor nt_d \rfloor + 1) dt \right)^2 \right].$$

Let $S_n(t) = n^{-d/2} \sum_{k \in A_n(t)} X_k$, then $M_n$ can be expressed as

$$M_n = \frac{1}{\hat{s}_n^2} \left[ \int_{[0,1]^d} \left( S_{n+1}(t) - \prod_{i=1}^d \frac{\lfloor nt_i \rfloor + 1}{n} S_n(1) \right)^2 dt \right.$$
$$\left. - \left( \int_{[0,1]^d} S_{n+1}(t) - \prod_{i=1}^d \frac{\lfloor nt_i \rfloor + 1}{n} S_n(1) dt \right)^2 \right].$$

Hence $\hat{s}_n^2 M_n$ is of the form $\Phi(S_n(.))$, where $\Phi$ is a continuous map.
As a consequence, under **SM2** and from the continuous mapping theorem,

$$\frac{\hat{s}_n^2}{\sigma^2} M_n \xrightarrow{\mathcal{L}} \int_{[0,1]^d} \left( B(t) - \left( \prod_{i=1}^d t_i \right) B(1) \right)^2 dt$$
$$- \left[ \int_{[0,1]^d} \left( B(t) - \left( \prod_{i=1}^d t_i \right) B(1) \right) dt \right]^2.$$

Besides, under **LM1**,

$$\hat{s}^2 \frac{n^d}{n^{2\gamma} L(n)^2} M_n \xrightarrow{\mathcal{L}} \int_{[0,1]^d} \left( Y(t) - \left( \prod_{i=1}^d t_i \right) Y(1) \right)^2 dt$$
$$- \left[ \int_{[0,1]^d} Y(t) - \left( \prod_{i=1}^d t_i \right) Y(1) dt \right]^2.$$

The proof is concluded thanks to the following lemma.

**Lemma 1.** *If* $\lim_{n\to\infty} q_n = \infty$ *and* $\lim_{n\to\infty} q_n/n = 0$, *then,*

(i) *Under* **H0***,*

$$\hat{s}_n^2 \xrightarrow{P} \sigma^2.$$

(ii) *Under* **H1***,*

$$\frac{n^d}{n^{2\gamma} L(n)^2} \hat{s}_n^2 \xrightarrow{P} 0.$$

*Proof of Lemma 1.*

The demonstration is an adaptation from Giraitis et al. (2003, [10]) to the random field framework.

Denote $k = (k_1, \ldots, k_d)$ and $|j| = (|j_1|, \ldots, |j_d|)$,

$$\tilde{r}(j) = \frac{1}{n^d} \sum_{k_1=1}^{n-|j_1|} \cdots \sum_{k_d=1}^{n-|j_d|} (X_k - \mu) \left( X_{k+|j|} - \mu \right),$$

where $\mu$ is the expectation of $X$. We split $\hat{s}_n^2$ as

$$\hat{s}_n^2 = \sum_{j \in B_{q-1}} \omega_{q,j} \tilde{r}(j) + \sum_{j \in B_{q-1}} \omega_{q,j} \left( \hat{r}(j) - \tilde{r}(j) \right) := u_n + v_n. \qquad (3.6)$$

We first show that, under **H0**, $E(|v_n|) \to 0$ when $n \to \infty$, while under **H1**, $E(|v_n|) = o(n^{2\gamma-d} L(n)^2)$ when $n \to \infty$.

Some computations lead to

$$\hat{r}(j) - \tilde{r}(j) = \prod_{i=1}^{d} \left( 1 - \frac{|j_i|}{n} \right) \left( \overline{X}_n - \mu \right)^2 - n^{-d} \left( \overline{X}_n - \mu \right)$$

$$\sum_{k_1=1}^{n-|j_1|} \cdots \sum_{k_d=1}^{n-|j_d|} \left( (X_k - \mu) + (X_{k+|j|} - \mu) \right).$$

From the Cauchy-Schwartz inequality

$$E(|v_n|) \leq \sum_{j \in B_{q-1}} E \left| \hat{r}(j) - \tilde{r}(j) \right|$$

$$\leq \sum_{j \in B_{q-1}} \left\{ E \left( \overline{X}_n - \mu \right)^2 + n^{-d} \sqrt{E \left( \overline{X}_n - \mu \right)^2} \left[ \sqrt{E \left( \sum_{k_1=1}^{n-|j_1|} \cdots \sum_{k_d=1}^{n-|j_d|} (X_k - \mu) \right)^2} \right. \right.$$

$$\left. \left. + \sqrt{E \left( \sum_{k_1=|j_1|}^{n} \cdots \sum_{k_d=|j_d|}^{n} (X_k - \mu) \right)^2} \right] \right\}.$$

Since $r$ is bounded, if $s_i < t_i$, for all $i = 1 \ldots d$,

$$E \left( \sum_{i_1=s_1}^{t_1} \cdots \sum_{i_d=s_d}^{t_d} (X_i - \mu) \right)^2 = \sum_{i_1, i_1'=s_1}^{t_1} \cdots \sum_{i_d, i_d'=s_d}^{t_d} r(i - i') \leq c \prod_{i=1}^{d} (t_i - s_i),$$

where $c$ is a positive constant. Therefore

$$E(|v_n|) \leq c \sum_{j \in B_{q-1}} \left( E\left(\overline{X}_n - \mu\right)^2 + 2n^{-d}\sqrt{E\left(\overline{X}_n - \mu\right)^2}\prod_{i=1}^{d}\sqrt{n - |j_i|}\right)$$

$$\leq c \sum_{j \in B_{q-1}} \left( E\left(\overline{X}_n - \mu\right)^2 + 2n^{-d/2}\sqrt{E\left(\overline{X}_n - \mu\right)^2}\right). \tag{3.7}$$

Under **SM1**, $E\left(\overline{X}_n - \mu\right)^2 \leq cn^{-d}$, while under **LM2**, $E\left(\overline{X}_n - \mu\right)^2 \leq cn^{2\gamma - 2d}L^2(n)$. Since $q/n$ goes to 0, it is then straightforward to conclude that, under **H0**, $E(|v_n|)$ vanishes and, under **H1**, $E(|v_n|) = o(n^{2\gamma - d}L^2(n))$.

Now, we turn to the asymptotic behavior of $u_n$.

Let us first show that under **H0**, $u_n$, defined by (3.6), converges in probability to $\sigma^2$ when $n$ goes to infinity. Notice that

$$E(u_n) = \sum_{j \in B_{q-1}} \omega_{q,j}\left(\prod_{i=1}^{d}\frac{n - |j_i|}{n}\right)r(j) \to \sigma^2.$$

According to (3.6), $\hat{s}_n^2$ converges in probability to $\sigma^2$ if $E\left(u_n - E(u_n)\right)^2$ converges to 0.

$$E\left(u_n - E(u_n)\right)^2$$

$$= E\left(\sum_{j \in B_{q-1}} \omega_{q,j}\left[\tilde{r}(j) - E\left(\tilde{r}(j)\right)\right]\right)^2$$

$$= \sum_{j,j' \in B_{q-1}^2} \omega_{q,j}\omega_{q,j'}\text{cov}(\tilde{r}(j), \tilde{r}(j'))$$

$$\leq \sum_{j,j' \in B_{q-1}^2} |\text{cov}(\tilde{r}(j), \tilde{r}(j'))|$$

$$\leq \frac{1}{n^{2d}} \sum_{j,j' \in B_{q-1}^2} \sum_{k,k' \in A_n^2} \left|\text{cov}\left((X_k - \mu)(X_{k+|j|} - \mu), (X_{k'} - \mu)(X_{k'+|j'|} - \mu)\right)\right|.$$

We split the sum above in two terms involving on one side the cumulants and on the other side the covariance function. Indeed

$$\text{cov}\left((X_k - \mu)(X_{k+|j|} - \mu), (X_{k'} - \mu)(X_{k'+|j'|} - \mu)\right)$$

$$= cum(X_k, X_{k+|j|}, X'_k, X_{k'+|j'|})$$

$$\quad + r(k - k')r(k' - k + |j'| - |j|) + r(k - k' - |j'|)r(k - k' + |j|).$$

First, thanks to **SM3**,

$$\frac{1}{n^{2d}} \sum_{j,j' \in B_{q-1}^2} \sum_{k,k' \in A_n^2} \left| cum(X_k, X_{k+|j|}, X'_k, X'_{k'+|j'|}) \right|$$

$$\leq \frac{1}{n^{2d}} \sum_{j \in B_{q-1}} \sum_{k \in A_n} \sum_{i,i' \in B_{2n}} \left| cum(X_0, X_{|j|}, X_i, X_{i'}) \right|$$

$$\leq \frac{1}{n^d} \sum_{j \in B_{q-1}} \sum_{i,i' \in B_{2n}} |c_4(|j|, i, i')|$$

$$\leq c \left( \frac{q}{n} \right)^d,$$

where $c$ is a positive constant.

Then,

$$\frac{1}{n^{2d}} \sum_{j,j' \in B_{q-1}^2} \sum_{k,k' \in A_n^2} |r(k-k')r(k'-k+|j'|-|j|) + r(k-k'-|j'|)r(k-k'+|j|)|$$

$$\leq \frac{1}{n^{2d}} \sum_{j \in B_{q-1}} \sum_{k \in A_n} \sum_{i,i' \in B_{2n}} 2 |r(i)r(i')|$$

$$\leq c \left( \frac{q}{n} \right)^d,$$

where $c$ is a positive constant.

Finally $E(u_n - E(u_n))^2$ converges to 0 if $q/n \to 0$ and this completes the proof of $(i)$.

For proving $(ii)$, recall that

$$E(u_n) = \sum_{j \in B_{q-1}} \omega_{q,j} \left( \prod_{i=1}^d \frac{n - |j_i|}{n} \right) r(j).$$

According to **LM2**, $E(u_n) \leq q^{-d} \sum_{j \in B_{q-1}} \left( \prod_{i=1}^d (q - |j_i|) \right) r(j) = O(q^{2\gamma-d} \times L^2(q))$. Therefore, from (3.6) and (3.7), $E(|\hat{s}_n^2|) = E(\hat{s}_n^2) = E(u_n) + E(v_n) = o(n^{2\gamma-d}L^2(n)) + O(q^{2\gamma-d}L^2(q))$ and $(ii)$ of Lemma 1 follows because $q/n$ vanishes when $n$ goes to infinity. $\qquad\square$

## 4. Simulations in dimension $d = 2$

The following simulations give an idea of the power of the test under different situations of strong dependence. First, one has to approach the asymptotic law of the test statistic $M_n$ under the null hypothesis. This is done in subsection 4.1. We then focus on the choice of $q$ when $n = 128$ and $n = 256$ to guarantee a proper size of the test. In subsection 4.3, the power of the test is assessed.

Several kinds of long memory random fields are simulated and submitted to the testing procedure. As these simulations are highly time consuming, we restrict ourselves to random fields of size $128 \times 128$ and $256 \times 256$. The simulations results reveal a close relation between the power and the kind of strong dependence encountered. The power depends on the strength of the long memory but also on its anisotropy.

### *4.1. Asymptotic law under* **H0**

The first step to implement the test consists in simulating the asymptotic law of $M_n$ under the null hypothesis. According to Proposition 1, this is the law of

$$\int_{[0,1]^2} \left( B(t) - \left( \prod_{i=1}^{2} t_i \right) B(1) \right)^2 dt - \left[ \int_{[0,1]^2} \left( B(t) - \left( \prod_{i=1}^{2} t_i \right) B(1) \right) dt \right]^2,$$
(4.1)

where $B$ is the Brownian Sheet on $[0,1]^2$. After some computations, (4.1) can be written

$$\int_{[0,1]^2} B(t_1, t_2)^2 dt_1 dt_2 - 2B(1,1) \int_{[0,1]^2} t_1 t_2 B(t_1, t_2) dt_1 dt_2$$

$$- \left( \int_{[0,1]^2} B(t_1, t_2) dt_1 dt_2 \right)^2 + \frac{B(1,1)}{2} \int_{[0,1]^2} B(t_1, t_2) dt_1 dt_2 + \frac{7}{144} B(1,1)^2.$$

To simulate a sample under this law, each integral above is approximated by a Riemann sum, for instance

$$\int_{[0,1]^2} t_1 t_2 B(t_1, t_2) dt_1 dt_2 \approx \frac{1}{n^2} \sum_{k_1=1}^{n} \sum_{k_2=1}^{n} \frac{k_1}{n} \frac{k_2}{n} B\left( \frac{k_1}{n}, \frac{k_2}{n} \right),$$

where a realization of $\left( B\left( \frac{k_1}{n}, \frac{k_2}{n} \right) \right)_{1 \le k_1, k_2 \le n}$ is given by

$$B\left( \frac{k_1}{n}, \frac{k_2}{n} \right) = \frac{1}{n} \sum_{j_1=1}^{k_1} \sum_{j_2=1}^{k_2} \varepsilon_{j_1, j_2}, \qquad \forall (k_1, k_2) \in \{1, \ldots, n\}^2,$$

with $(\varepsilon_j)_{j \in \mathbb{Z}^2}$ a Gaussian white noise.

For $n = 7000$, 10000 realizations of the law (4.1) have been computed. The histogram of the sample is shown in Figure 1.

From the simulated sample, an estimated mean of 0.0897 is obtained, the empirical variance is 0.0018 and the empirical quantiles of order 90% and 95% are respectively 0.1448 and 0.1692.

### *4.2. Choice of* $q$

The sample size $n$ being fixed, one has to choose the value of $q$ involved in definition (3.3) of $M_n$. This choice is not easy. We decide to choose it so that the size of the test is optimized. This reduces to check that the p-values under **H0** are uniformly distributed on $[0, 1]$.

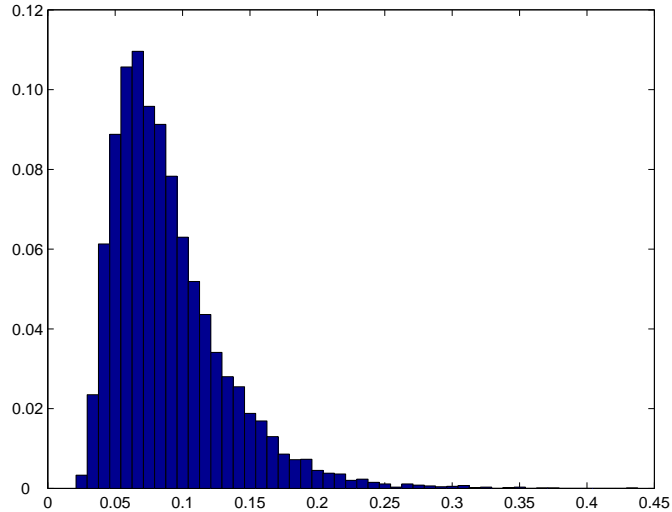

FIG 1. *Estimated density function of the limiting law (4.1).*

Small values of $q$ lead to an increase of the probability of rejecting **H0**. So, in order to maximize the power of the test, we have to find the smallest $q$ such that the size is correct.

To achieve this choice, we compute the test for different autoregressive fields $X_{k_1, k_2}$ defined by

$$(1 - aL_1)(1 - aL_2)X_{k_1, k_2} = \varepsilon_{k_1, k_2}, \tag{4.2}$$

where $\varepsilon$ is a Gaussian white noise, where $0 < a < 1$ and $L_i$ represents the lag operator on the i-*th* index. The simulation of such an autoregressive field is done by filtering a white noise as in [8].

Clearly, the size of the test will increase with $a$. As a consequence, in the simulations below, we choose to quote only two situations: $a = 0.5$ and $a = 0.8$. The last case corresponds to a memory close to strong dependence and a worse size in this case might be acceptable. Indeed, it will be impossible to find $q$ such that, uniformly on $a$, the size of the test is strictly the one expected.

### 4.2.1. Case $n = 128$

Figure 2 represents the empirical distribution of the p-values, computed on 1000 realizations of $M_n$, in the case $n = 128$, when $q = 28$, $q = 30$ and $q = 32$. The line with crosses stands for $a = 0.8$ and the line with circles for $a = 0.5$. The diagonal line is added for sake of comparison. The representations are zoomed in on $[0, 0.1]$.

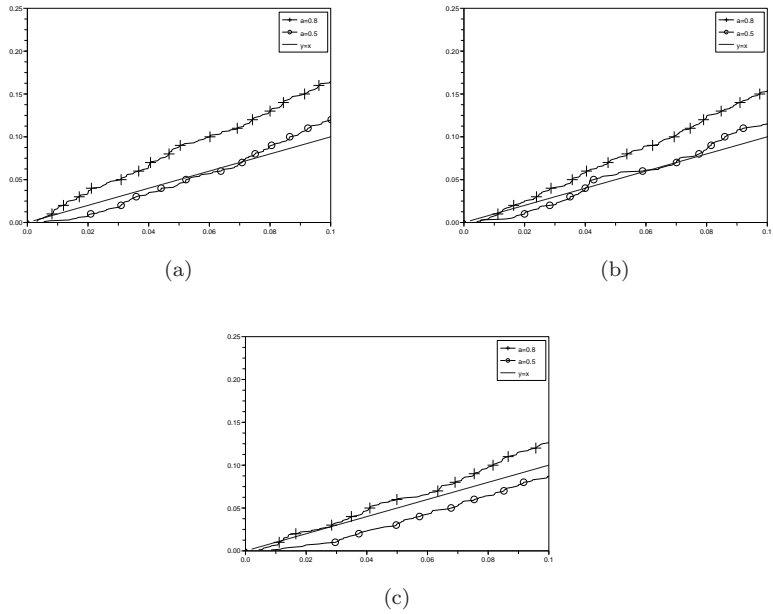

FIG 2. *Cumulative distribution function of the p-values on* $[0, 0.1]$ *for model* (4.2) *with* $a = 0.8$ *(crosses) and* $a = 0.5$ *(circles) when* $n = 128$ *and* $q = 28$ *(a),* $q = 30$ *(b),* $q = 32$ *(c).*

The value $q = 30$ is chosen. For this value, the size associated with $a = 0.8$ is larger than expected (for instance 15% instead of 10%). This error is acceptable since the dependence of an autoregressive field (4.2) with $a = 0.8$ is close to long memory. The choice of a larger $q$ reduces this error but on the other hand, this may create a bias for smaller values of $a$ (as seen in Figure 2 for $q = 32$ and $a = 0.5$), and consequently the test becomes less powerful.

### 4.2.2. Case $n = 256$

Figure 3 represents the empirical distribution of the p-values in the case $n = 256$ when $q = 35$, $q = 40$ and $q = 45$. The line with crosses stands for $a = 0.8$ and the line with circles for $a = 0.5$.

The value $q = 40$ is chosen for the same reasons as before: The error for $a = 0.8$ seems acceptable in comparison with the loss of power that the choice of a larger $q$ would yield.

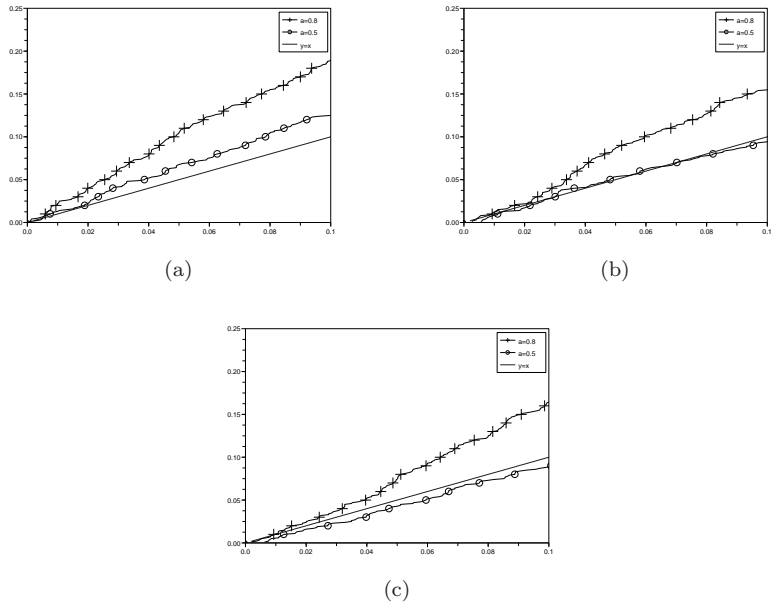

Fig 3. *Cumulative distribution function of the p-values on $[0, 0.1]$ for model (4.2) with $a = 0.8$ (crosses) and $a = 0.5$ (circles) when $n = 256$ and $q = 35$ (a), $q = 40$ (b), $q = 45$ (c).*

## *4.3. Power under different alternatives*

We implement the test on different Gaussian long memory random fields. The power is assessed according to the type of memory. First, the case when the long memory is of tensorial product type is studied: That is when the spectral density $f(x_1, x_2)$ of the field, defined on $[-\pi, \pi]^2$, is equivalent at 0 to $|x_1|^{\alpha_1} |x_2|^{\alpha_2}$, $-1 < \alpha_1 < 0$, $-1 < \alpha_2 < 0$. The range of $\alpha_1$ and $\alpha_2$ guarantees the integrability of $f$. Notice that this spectral density does not follow the property of Definition 1. Then, we report the case when the long memory is isotropic: That is when the spectral density is equivalent at 0 to $(x_1^2 + x_2^2)^{\alpha/2}$, where, for integrability reasons, $-2 < \alpha < 0$. In the last simulations, we focus on the effect of anisotropy on the power. For this purpose, we simulate Gaussian fields whose spectral densities are equivalent at 0 to $|x_1 + kx_2|^\alpha$, $-1 < \alpha < 0$, $k \in \mathbb{Z}$. This is an extreme case of anisotropic field since $f$ exhibits only one line of singularity. For all these examples, assumption **H1** holds with $\gamma = 1 - (\alpha_1 + \alpha_2)/2$ for the product-type case, and $\gamma = 1 - \alpha/2$ for the two other examples (see [14] or [15]).

All of the above fields are simulated thanks to the spectral method (cf. [13]) which consists in the following algorithm:

1. Generate N independent random variables $(Z_1^{(1)}, Z_2^{(1)}), \dots, (Z_1^{(N)}, Z_2^{(N)})$ on $[-\pi, \pi]^2$ according to the spectral measure $\mu$ (viewed, up to a normalization, as a probability distribution);

2. Generate N independent random variables $U_1, \ldots, U_N$ uniformly on $[0,1]$;
3. Compute for all $(i,j)$, $X_{i,j} = \frac{\sqrt{2}}{N} \sum_{k=1}^{N} \cos(Z_1^{(k)} i + Z_2^{(k)} j + 2\pi U_k)$.

According to the central limit theorem, the resulting field $X$ is Gaussian with spectral measure $\mu$ when $N$ is large. In practice, the value $N = 5000$ is fixed.

### 4.3.1. Tensorial product type long memory

The power of the test is assessed for Gaussian fields with a spectral density equivalent at 0 to $|x_1|^{\alpha_1}|x_2|^{\alpha_2}$, $-1 < \alpha_1 < 0$, $-1 < \alpha_2 < 0$. Figure 4 shows a simulation of such fields on a $256 \times 256$ grid using the spectral method, when, from left to right, $\alpha_1 = \alpha_2 = -0.25$, $\alpha_1 = \alpha_2 = -0.5$ and $\alpha_1 = \alpha_2 = -0.75$. These cases correspond respectively in (2.5) to $\gamma = 1.25$, $\gamma = 1.5$ and $\gamma = 1.75$.

The empirical c.d.f. of the p-values of the test is represented on Figure 5. Each curve is computed on 500 simulated fields. On the left, the simulations correspond to fields of size $128 \times 128$, while on the right $n = 256$. The parameter $q$ has been chosen according to subsection 4.2, that is $q = 30$ when $n = 128$ and $q = 40$ when $n = 256$. In each case, three curves are represented, corresponding to $\alpha_1 = \alpha_2 = -0.25$ (crosses), $\alpha_1 = \alpha_2 = -0.5$ (circles) and $\alpha_1 = \alpha_2 = -0.75$ (triangles). We observe logically that the power of the test increases with the strength of the memory (quantified by $\gamma$) and with the size of the sample.

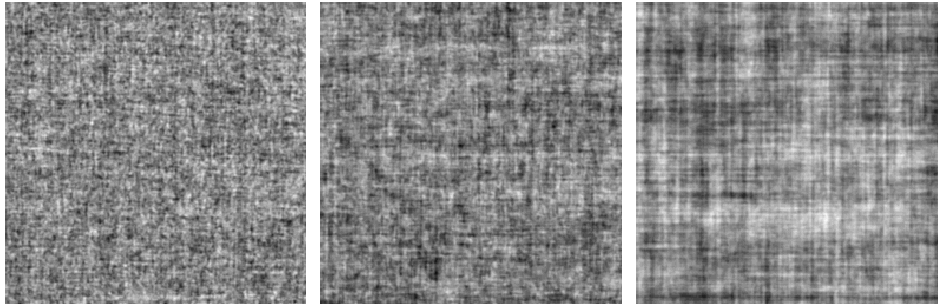

FIG 4. *Gaussian fields with product-type long memory where $\gamma = 1.25$ (left), $\gamma = 1.5$ (middle) and $\gamma = 1.75$ (right).*

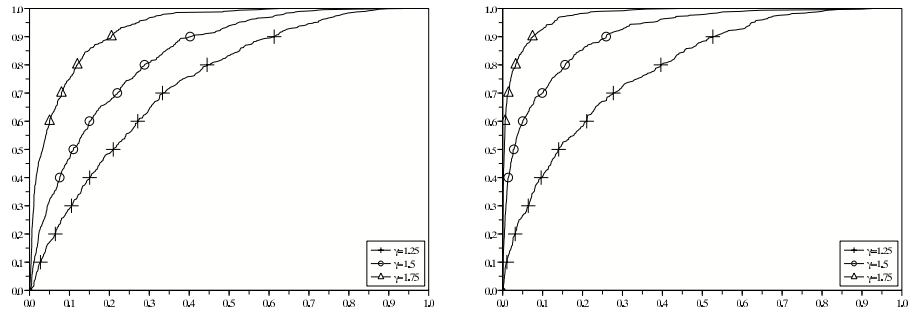

FIG 5. *C.d.f. of the p-values when $\gamma = 1.25$ (crosses), $\gamma = 1.5$ (circles) and $\gamma = 1.75$ (triangles), where $n = 128$ (left) and $n = 256$ (right).*

### 4.3.2. Isotropic long memory

Now, the power of the test is assessed on isotropic long memory fields. Their spectral density is equivalent at 0 to $(x^2 + y^2)^{\alpha/2}$, $-2 < \alpha < 0$. A simulation of such fields is presented on Figure 6 when $\alpha = -0.5$, $\alpha = -1$ and $\alpha = -1.5$, which correspond to $\gamma = 1.25$, $\gamma = 1.5$ and $\gamma = 1.75$ in (2.5).

Figure 7 represents the same c.d.f as in Figure 5 but for the isotropic case when $\alpha = -0.5$ (crosses), $\alpha = -1$ (circles) and $\alpha = -1.5$ (triangles). These choices correspond to the same strengths of long memory (i.e. the same $\gamma$) as in Figure 5. The power follows a similar behavior: It increases with the strength of the memory and with the sample size. However, the power for the isotropic case is smaller than in the product-type setting studied before. This is due to the form of the statistic $M_n$ (see Section 3.1): The empirical variance of $S_{n,j}^*$ is computed using quadrants $A_j$'s. This suits better product-type fields than isotropic ones for detecting long memory. Indeed, in the product-type setting, this empirical variance will tend to return higher values. This sensitivity to anisotropy is studied further in the last subsection.

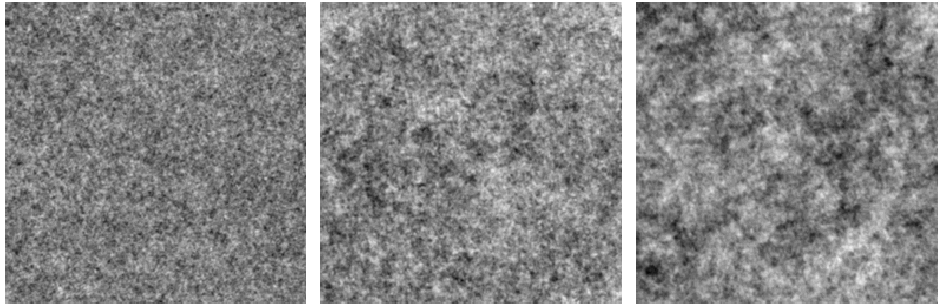

FIG 6. *Gaussian fields with isotropic long memory when $\gamma = 1.25$ (left), $\gamma = 1.5$ (middle) and $\gamma = 1.75$ (right).*

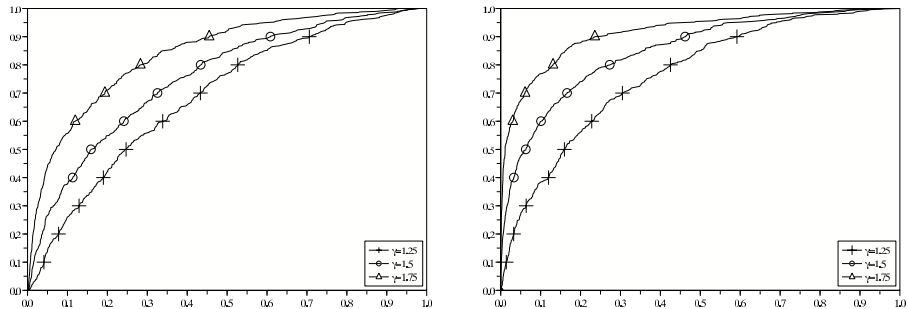

FIG 7. *C.d.f. of the p-values when $\gamma = 1.25$ (crosses), $\gamma = 1.5$ (circles) and $\gamma = 1.75$ (triangles), where $n = 128$ (right) and $n = 256$ (left).*

### *4.3.3. The effect of anisotropy on the power*

As seen before, the power for the isotropic long memory case is smaller than for the product-type setting. To assess properly the sensitivity to anisotropy, we focus on one-direction long memory fields in the sense that their spectral density behaves at zero as $|x_1 + kx_2|^\alpha$, $k \in \mathbb{Z}$. Note that this form of the spectral density demands $-1 < \alpha < 0$ to be integrable. This yields the restriction $\gamma < 1.5$ in (2.5). Figure 8 shows a simulation on a $256 \times 256$ grid when $k = -1$ and $k = 0$ and when $\alpha = -0.5$, corresponding to $\gamma = 1.25$ in (2.5). As we can see on Figure 8, the worst case, in terms of computation using quadrants, should be $k = -1$ while the most suitable situation should be $k = 0$.

Figure 9 shows the power of the test when $k = -1$ (solid line) and $k = 0$ (dotted line) for $n = 128$ (left) and $n = 256$ (right). This representation confirms the effect of anisotropy. Therefore, from a practical point of view, it seems better to study first the isotropy of the dependence before testing the presence of long memory, in order to rotate the image sample if necessary.

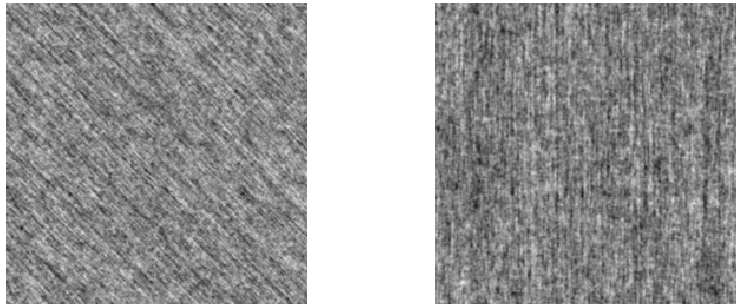

FIG 8. *Gaussian fields with one-direction long memory ($\gamma = 1.25$) when $k = -1$ (left) and $k = 0$ (right).*

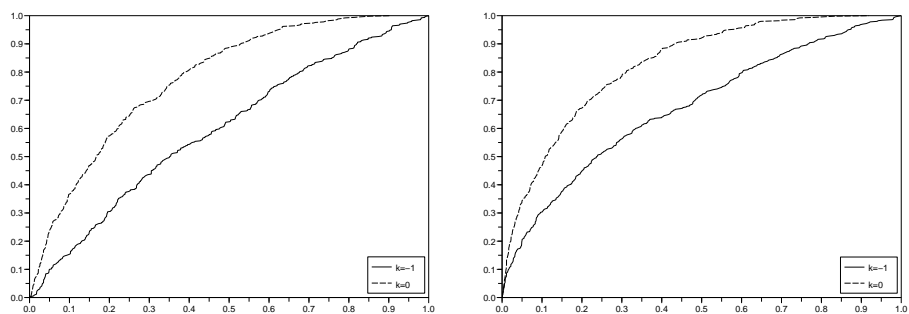

Fig 9. *C.d.f. of the p-values when $k = -1$ (solid line) and $k = 0$ (dotted line) for $\gamma = 1.25$ and where $n = 128$ (left) and $n = 256$ (right).*